\documentclass[preprint]{elsarticle}

\usepackage{amssymb}
\usepackage{amsmath}
\usepackage{amsthm}
\usepackage{enumerate}
\usepackage{graphicx}
\usepackage{amsfonts}

\journal{...}

\newtheorem{thm}{Theorem}

\newtheorem{tb}{Table}

\theoremstyle{definition}

\newtheorem{remark}{Remark}

\begin{document}

\begin{frontmatter}

\title{Generalized Alomari functionals}

 \author[1]{ Ana-Maria Acu}
  \author[2]{Heiner Gonska}
 \address[1]{Lucian Blaga University of Sibiu, Department of Mathematics and Informatics, Str. Dr. I. Ratiu, No.5-7, RO-550012  Sibiu, Romania,
 e-mail: acuana77@yahoo.com}
\address[2]{University of Duisburg-Essen, Faculty of Mathematics, Forsthausweg 2, 47057 Duisburg, Germany, e-mail: heiner.gonska@uni-due.de }

\begin{abstract}
We consider a generalization form of certain integral inequalities
given by Guessab, Schmeisser and Alomari. The trapezoidal, mid
point, Simpson, Newton-Simpson rules are obtained as special cases.
Also,  inequalities for the generalized Alomari functional in term
of the $n$-th order modulus, $n=\overline{1,4}$, are given and
applied to some known quadrature rules.

\end{abstract}


\end{frontmatter}

\section{Introduction}
In the present note we will consider a certain quadrature functional defined for functions in $(C[a,b],\|\cdot\|_{\infty})$, the space of continuous functions defined on the compact interval $[a,b]$, $a<b$, equipped with the sup norm $\|\cdot\|_{\infty}$.

In 1938, Ostrowski \cite{Ostrowski} published one of the classical
inequalities dealing with the most primitive form of a quadrature rule.
The  Ostrowski inequality gives an  approximation of an integral by a
single value of the function.
\begin{thm} Let $f$ be a differentiable function on $(a,b)$, and let
$|f^{\prime}(t)|\leq M$ for $t\in (a,b)$. Then, for each $x\in
(a,b)$,
\begin{equation}\label{e1}
\left|\displaystyle\frac{1}{b-a}\int_{a}^bf(t)dt-f(x)\right|\leq
\left[\displaystyle\frac{1}{4}+\left(\frac{x-\frac{1}{2}(a+b)}{b-a}\right)^2\right](b-a)M.
\end{equation}
The constant $\displaystyle\frac{1}{4}$ is the best possible in the
sense that it cannot be replaced by a smaller constant.
\end{thm}
Ostrowski inequalities have attracted the attention of many
mathematicians. The reader should consult the monographies  of D.S.
Mitrinovi\'{c} et al. \cite{Mitrinovic}, G. Anastassiou
\cite{Anastassiou} and the recent book by P.Cerone et al. \cite{Cerone}.

Let us consider the  space of Lipschitz function of order
$\alpha\in(0,1]$ and the constant $M>0$, which is defined as follows
$$Lip_{M}(\alpha):=\left\{f\in C[a,b]:\,|f(x)-f(y)|\leq M|x-y|^{\alpha}, \textrm{for all } x,y\in[a,b]\right\}.$$
The class  $Lip_{M}(1)$ is simply denote by $Lip_{M}$.

 For the
trapezoidal rule and mid point rule, S.S. Dragomir at al.
\cite{Dragomir} obtained the following inequalities:
\begin{thm}\label{t2}\cite{Dragomir} Let $f$ be a function defined on an interval $[a,b]$ and
belonging to $Lip_{M}$. Then
\begin{equation}\label{a} \left|\displaystyle\frac{1}{b-a}\int_{a}^b f(t)dt-f\left(\frac{a+b}{2}\right)\right|\leq
\displaystyle\frac{M}{4}(b-a), \end{equation}
and
\begin{equation}\label{b} \left|\displaystyle\frac{1}{b-a}\int_{a}^b f(t)dt-\frac{f(a)+f(b)}{2}\right|\leq \displaystyle\frac{M}{3}(b-a).\end{equation}
\end{thm}
In 2002, A. Guessab and  G. Schmeisser \cite{Guessab} introduced a
new integral inequality which can be interpreted as an analogue of
Ostrowski's result.
\begin{thm}\cite{Guessab} Let $f$ be a function defined on $[a,b]$
and belonging to $Lip_{M}(\alpha)$ with $\alpha\in(0,1]$. Then, for
each $x\in\left[a,\frac{1}{2}(a+b)\right]$, it follows
\begin{equation}\label{e2} \left|\displaystyle\frac{1}{b-a}\int_a^b
f(t)dt-\frac{f(x)+f(a+b-x)}{2}\right|\leq\displaystyle\frac{M}{b-a}\cdot\frac{(2x-2a)^{\alpha+1}+(a+b-2x)^{\alpha+1}}{2^\alpha(\alpha+1)}.
\end{equation}
This inequality is sharp for each admissible $x$. Equality  is
attained if and only if $f=\pm Mf_*+c$ with $c\in \mathbb{R}$ and
$$ f_*(t):=\left\{\begin{array}{lll} (x-t)^\alpha,&\textrm{for}&a\leq t\leq x,\\
\vspace{-0.3cm}\\
(t-x)^\alpha& \textrm{for}&x\leq t\leq \frac{1}{2}(a+b),\\
 \vspace{-0.3cm}\\
 f_*(a+b-t)&\textrm{for}& \frac{1}{2}(a+b)\leq t\leq b.\end{array}\right.
$$
\end{thm}
\begin{remark} Setting $\alpha=1$ in the inequality (\ref{e2}),
the following result  is derived:
 \begin{equation}\label{e3} \left|\displaystyle
\frac{1}{b-a}\int_a^b f(t) dt-\frac{f(x)+f(a+b-x)}{2}\right|\leq
\left[\frac{1}{8}+2\left(\frac{x-\frac{3a+b}{4}}{b-a}\right)^2\right](b-a)M,
\end{equation} where $f\in Lip_M$ and
$x\in\left[a,\displaystyle\frac{a+b}{2}\right]$. The constant
$\displaystyle\frac{1}{8}$ is the best possible in the sense that it
cannot be replaced by a smaller constant.
\end{remark}
\begin{remark} Setting $x=\displaystyle\frac{1}{2}(a+b)$ in the inequality (\ref{e3}) we recover the estimates (\ref{a}) of Theorem \ref{t2}. For $x=a$ the estimation (\ref{b}) is improved, obtaining the following estimation of remainder term in trapezoidal rule
$$ \left|\displaystyle\frac{1}{b-a}\int_{a}^b f(t)dt-\frac{f(a)+f(b)}{2}\right|\leq \displaystyle\frac{M}{4}(b-a).$$
\end{remark}
Also, in \cite{Alomari} Alomari obtained the following
\begin{thm}\label{t4}\cite{Alomari} Let $f$ be a differentiable function on $(a,b)$, and let
$|f^{\prime}(t)|\leq M$ for $t\in (a,b)$. Then, the inequality
\begin{align}
&\left|\displaystyle\frac{1}{(b-a)}\int_a^b
f(t)dt-\left[\lambda\frac{f(a)+f(b)}{2}+(1-\lambda)\frac{f(x)+f(a+b-x)}{2}\right]\right|\nonumber\\
&\leq
(b-a)\left[\displaystyle\frac{1}{8}\left(2\lambda^2+(1-\lambda)^2\right)
+2\left(\frac{x-\frac{(3-\lambda)a+(1+\lambda)b}{4}}{b-a}\right)^2\right]\|f^{\prime}\|_{\infty}\label{e4}
\end{align}
holds for all $\displaystyle\lambda\in[0,1]$ and $\displaystyle a+\lambda\frac{b-a}{2}\leq x\leq \frac{a+b}{2}$.
\end{thm}
The functional in Theorem \ref{t4} is the Alomari functional.
\begin{remark} For $\lambda=0$ it is recovered the Guessab and
Schmeisser' s result (\ref{e3}). Also, using the inequality
(\ref{e4}) are obtained estimates for remainder term of the
midpoint, trapezoid and Simpson formulas.
\end{remark}

\section{Estimates for the generalized Alomari functional}
The previous results motivated us to study, for all real numbers
$x\in [a,b]$ and all real parameters $ \lambda,\mu\in[0,1] $, the
remainder term of the more general Alomari quadrature formula as given between the accolades:
\begin{equation}\label{e5}
E(f;\lambda,\mu;x):=\displaystyle\frac{1}{b-a}\int_a^b f(t)
dt-\left\{\lambda\frac{f(a)+f(b)}{2}+(1-\lambda)\cdot\left[\mu
f(x)+(1-\mu)f(a+b-x)\right]\right\}.
\end{equation}

We observe that (\ref{e5}) defines a family of quadrature formulas
which contains the trapezoidal, mid point, Simpson rules. Also, it
is easy to see that for $\lambda=0,\mu=\displaystyle\frac{1}{2}$ we
recover Guessab and Schmeisser's functional in (\ref{e3}) and for
$\mu=\displaystyle\frac{1}{2}$ we have the one of Alomari in
(\ref{e4}).

\subsection{Estimates of the functional $E(f;\lambda,\mu;x)$ for continuously differentiable functions}
Let
$$ F(f)=\displaystyle\int_a^bf(x)dx-\displaystyle\sum_{k=1}^{n} A_kf(x_k) , \textrm{ where } A_k\in \mathbb{R} \textrm{ and } x_k\in[a,b]$$
be a linear functional which vanishes on $\Pi_k$ (polynomials of degree $\leq k$).

Using Taylor's formula, for a sufficiently smooth function $f$ defined on $[a,b]$, it follows
\begin{equation}\label{XX}
f=\displaystyle\sum_{j=0}^kf^{(j)}(a)\displaystyle\frac{(\cdot-a)_{+}^j}{j!}+\int_a^bf^{(k+1)}(t)\displaystyle\frac{(\cdot-t)_{+}^k}{k!}dt,
\end{equation}
where $(\cdot)_{+}^k$ is the truncaded power function. Applying $F$ to (\ref{XX}) gives
\begin{align*}
F(f)&= F\left(\displaystyle\int_a^bf^{(k+1)}(t)\frac{(\cdot-t)_{+}^k}{k!}dt\right)\\
&=\displaystyle\int_a^b\int_a^b f^{(k+1)}(t)\frac{(x-t)_{+}^k}{k!}dtdx-\sum_{k=1}^nA_k\int_a^bf^{(k+1)}(t)\frac{(x_k-t)_{+}^k}{k!}dt\\
&=\displaystyle\int_a^bf^{(k+1)}(t)\left[\int_a^b\frac{(x-t)_{+}^k}{k!}dx-\sum_{k=1}^nA_k\frac{(x_k-t)_{+}^k}{k!}\right]dt=\displaystyle\frac{1}{k!}\int_a^bf^{(k+1)}(t)F\left((\cdot-t)_{+}^k\right)dt.
\end{align*}
If we note $K(t)=\displaystyle\frac{1}{k!}F\left((\cdot-t)_{+}^k\right)$ we have the following representation of the functional $F$:
$$  F(f)=\displaystyle\int_a^bf^{(k+1)}(t)K(t)dt.$$
The function $K$ is known as Peano's  kernel (see \cite{Peano}).

\noindent For all $\lambda,\mu\in [0,1]$ and $a\leq x\leq
\displaystyle\frac{a+b}{2}$, define the Peano's kernel
$$K(x,t;\lambda,\mu):=\left\{\begin{array}{l} t-\nu_1(\lambda),\, t\in[a,x],\\
\vspace{-0.3cm}\\
t-\nu_2(\lambda,\mu),\,t\in(x,a+b-x),\\
\vspace{-0.3cm}\\
t-\nu_3(\lambda),\,t\in[a+b-x,b],\end{array}\right.  $$
where
\begin{align*} &\nu_1(\lambda):=a+\displaystyle\frac{\lambda}{2}(b-a),
\nu_2(\lambda,\mu):=\displaystyle
a\!+\!(b\!-\!a)\left(\frac{\lambda}{2}\!+\!\mu(1\!-\!\lambda)\right),\nu_3(\lambda):=b\!-\!\displaystyle\frac{\lambda}{2}(b\!-\!a)\!=\!a\!+\!b\!-\!\nu_1(\lambda).\end{align*}
Integrating by parts, we get
\begin{align*}\displaystyle\int_a^b K(x,t;\lambda,\mu) f^{\prime}(t)dt\!&\!=\!(b\!-\!a)\left\{\lambda\displaystyle\frac{f(a)\!+\!f(b)}{2}
\!+\!(1\!-\!\lambda)\left[\mu f(x)\!+\!(1\!-\!\mu)f(a\!+\!b\!-\!x)\right]\!-\!\frac{1}{b\!-\!a}\int_a^bf(t) dt\right\}\\
&=-(b-a)E(f;\lambda,\mu;x).
\end{align*}
Therefore,
\begin{equation}\label{e7}
E(f;\lambda,\mu;x)=-\displaystyle\frac{1}{b-a}\int_a^bK(x,t;\lambda,\mu) f^{\prime}(t)dt.
\end{equation}
\begin{thm}\label{t5}Let $f:[a,b]\to \mathbb{R}$ be differentiable on the
interval $(a,b)$, with the first derivative bounded on $(a,b)$, i.e.
$\|f^{\prime}\|:=\displaystyle\sup_{t\in[a,b]}|f^{\prime}(t)|<\infty
$.Then, for all $\lambda,\mu\in[0,1]$ and
$x\in\left[a,\frac{a+b}{2}\right]$, we have
\begin{equation}\label{e6}
|E(f;\lambda,\mu;x)|\leq\displaystyle\frac{C(x,\lambda,\mu)}{2(b-a)}\cdot\|
f^{\prime}\|_{\infty},\,\textrm{where}
\end{equation}
$$C(x,\lambda,\mu):=2\lambda(x-a)(b-a)-2(x-a)^2+(x-\nu_2(\lambda,\mu)^2+(a+b-x-\nu_2(\lambda,\mu)^2,
\textrm{ for } x\in\left[a,a+\lambda\frac{b-a}{2}\right]$$ and
 $$C(x,\lambda,\mu):=(x-a)^2+\left[\lambda(b-a)-(x-a)\right]^2+\tilde{C}(x,\lambda,\mu), \textrm{ for } x\in\left[a+\lambda\frac{b-a}{2},\frac{a+b}{2}\right],$$
 with
 $$\tilde{C}(x,\lambda,\mu):=\left\{\begin{array}{ll}
(x-\nu_2(\lambda,\mu))^2+(a+b-x-\nu_2(\lambda,\mu))^2, &\textrm{if }
\mu\in\left[0,\frac{1}{2}\right], x\leq
\nu_2(\lambda,\mu),\\
\vspace{-0.3cm}\\
-(x-\nu_2(\lambda,\mu))^2+(a+b-x-\nu_2(\lambda,\mu))^2, &\textrm{if
} \mu\in\left[0,\frac{1}{2}\right], x>
\nu_2(\lambda,\mu),\\
\vspace{-0.3cm}\\
(x-\nu_2(\lambda,\mu))^2+(a+b-x-\nu_2(\lambda,\mu))^2, &\textrm{if }
\mu\in\left(\frac{1}{2},1\right], x\leq
a+b-\nu_2(\lambda,\mu),\\
\vspace{-0.3cm}\\
(x-\nu_2(\lambda,\mu))^2-(a+b-x-\nu_2(\lambda,\mu))^2, &\textrm{if }
\mu\in\left(\frac{1}{2},1\right], x>
a+b-\nu_2(\lambda,\mu).\\
\end{array}
 \right.  $$
\end{thm}

$   $

\begin{proof} Using the equality (\ref{e7}) we obtain the following estimation of the remainder term
$$ |E(f;\lambda,\mu;x)|\leq\displaystyle\frac{1}{b-a}\|f^{\prime}\|_{\infty}\int_a^b|K(x,t;\lambda,\mu)|dt. $$
Denote
$
I_1:=\displaystyle\int_a^x|t-\nu_1(\lambda)|dt,\,
I_2:=\displaystyle\int_x^{a+b-x}|t-\nu_2(\lambda,\mu)|dt,\,
I_3:=\displaystyle\int_{a+b-x}^b|t-\nu_3(\lambda)| dt.
$

For $x\in\left[a,a+\lambda\displaystyle\frac{b-a}{2}\right)$, it follows
\begin{align*}
I_1&=\displaystyle\int_a^x|t-\nu_1(\lambda)|dt=-\int_a^x[t-\nu_1(\lambda)]dt=\frac{1}{2}\left\{\lambda(x-a)(b-a)-(x-a)^2\right\},\\
I_2&=\displaystyle\int_x^{a+b-x}|t-\nu_2(\lambda,\mu)|dt=-\int_x^{\nu_2(\lambda,\mu)}\left[t-\nu_2(\lambda,\mu)\right]dt+\int_{\nu_2(\lambda,\mu)}^{a+b-x}
\left[t-\nu_2(\lambda,\mu)\right]dt\\
&=\displaystyle\frac{1}{2}\left[x-\nu_2(\lambda,\mu)\right]^2+\frac{1}{2}\left[a+b-x-\nu_2(\lambda,\mu)\right]^2,\\
I_3&=\displaystyle\int_{a+b-x}^b|t-\nu_3(\lambda)|dt=-\frac{1}{2}\left\{(x-a)^2+\lambda(a-x)(b-a)\right\}.
\end{align*}

For $x\in\left[a+\frac{\lambda}{2}(b-a),\frac{a+b}{2}\right]$, by direct calculations, we get
\begin{align*}
I_1&=-\int_a^{\nu_1(\lambda)}\left[t-\nu_1(\lambda)\right]dt+\int_{\nu_1(\lambda)}^x\left[t-\nu_1(\lambda)\right] dt
=\frac{1}{4}(x-a)^2+\left(\lambda\frac{b-a}{2}-\frac{x-a}{2}\right)^2,\\
I_3&=-\int_{a+b-x}^{\nu_3(\lambda)}\left[t-\nu_3(\lambda)\right]dt+\int_{\nu_3(\lambda)}^b\left[t-\nu_3(\lambda)\right] dt
=\frac{1}{4}(x-a)^2+\left(\lambda\frac{b-a}{2}-\frac{x-a}{2}\right)^2.
\end{align*}
In order to calculate $I_2$ we consider two cases $\displaystyle\mu\in\left[0,\frac{1}{2}\right]$ and $\displaystyle\mu\in\left(\frac{1}{2},1\right]$.

\noindent{\bf Case 1.} If $\displaystyle\mu\in\left[0,\frac{1}{2}\right]$, then $\displaystyle a+\frac{\lambda}{2}(b-a)\leq\nu_2(\lambda,\mu)\leq\frac{a+b}{2}$,
 and we have
 $$ I_2=\left\{\begin{array}{l}\displaystyle\frac{1}{2}\left(x-\nu_2(\lambda,\mu)\right)^2+\frac{1}{2}\left(a+b-x-\nu_2(\lambda,\mu)\right)^2,
 \textrm{ if } x\leq \nu_2(\lambda,\mu),\\
 \vspace{-0.3cm}\\
 \displaystyle\frac{1}{2}\left(a+b-x-\nu_2(\lambda,\mu)\right)^2-\frac{1}{2}\left(x-\nu_2(\lambda,\mu)\right)^2,
 \textrm{ if } x> \nu_2(\lambda,\mu).\\
 \end{array}\right. $$

\noindent{\bf Case 2.} If $\displaystyle\mu\in\left(\frac{1}{2},1\right]$, then $\displaystyle\frac{a+b}{2}<\nu_2(\lambda,\mu)\leq b-\frac{\lambda}{2}(b-a)$,
 and we have
 $$ I_2=\left\{\begin{array}{l}\displaystyle\frac{1}{2}\left(x-\nu_2(\lambda,\mu)\right)^2+\frac{1}{2}\left(a+b-x-\nu_2(\lambda,\mu)\right)^2,
 \textrm{ if } x\leq a+b-\nu_2(\lambda,\mu),\\
 \vspace{-0.3cm}\\
 \displaystyle-\frac{1}{2}\left(a+b-x-\nu_2(\lambda,\mu)\right)^2+\frac{1}{2}\left(x-\nu_2(\lambda,\mu)\right)^2,
 \textrm{ if } x> a+b-\nu_2(\lambda,\mu).\\
 \end{array}\right. $$
 Using the values of $I_k, k=1,2,3$, the inequality (\ref{e6}) is proved.
\end{proof}
\begin{remark}For $x\in\left(\displaystyle\frac{a+b}{2},b\right]$ we obtain the following representation of the functional $E$:
$$
E(f;\lambda,\mu;x)=-\displaystyle\frac{1}{b-a}\int_a^b\tilde{K}(x,t;\lambda,\mu) f^{\prime}(t)dt,
$$
where
$$\tilde{K}(x,t;\lambda,\mu):=\left\{\begin{array}{l} t-\nu_1(\lambda),\, t\in[a,a+b-x],\\
\vspace{-0.3cm}\\
t-\left[a+b-\nu_2(\lambda,\mu)\right],\,t\in(a+b-x,x),\\
\vspace{-0.3cm}\\
t-\nu_3(\lambda),\,t\in[x,b].\end{array}\right.  $$
Therefore, for $x\in\left(\displaystyle\frac{a+b}{2},b\right]$ can be obtained the similar estimation of $E(f;\lambda,\mu;x)$ as in\linebreak Theorem \ref{t5}.
\end{remark}

\subsection{Special cases}

We will next consider in which estimates in terms of the norms of higher order derivatives are possible.
 Each special case below is such that the quadrature rule reaches its maximal degree of exactness as indicated by the order of the derivative.


In order to get all these special cases we will calculate the maximum degree of exactness of the functional $E$, namely we determine the integer and positive number $n$ such that
$E(e_i;\lambda,\mu;x)=0$, for $i=\overline{0,n}$  and $E(e_{n+1};\lambda,\mu;x)\ne 0  $, where $e_i(x)=x^i$. All these cases are considered in the below table:
\begin{tb}\label{table1}
{\it The maximum degree n of exactness of the functional $E$}

$  $

${
\begin{tabular}{|l|l|}\hline
\textrm{case} &$n$\\ \hline
 $\displaystyle\lambda\ne 1, \mu\ne\frac{1}{2},x\ne\frac{a+b}{2}$& $0$\\
 \hline
 $\displaystyle\lambda=1$&$1$\\
 \hline
 $\displaystyle x=\frac{a+b}{2},\lambda\ne\frac{1}{3}$ &$1$\\
  \hline
 $\displaystyle\mu=\frac{1}{2},\lambda\in\left(\frac{1}{3},1\right)$&$1$\\
 \hline
 $\mu=\frac{1}{2},\lambda\in\left[0,\frac{1}{3}\right),x\ne\displaystyle\frac{a+b}{2}-\frac{b-a}{2}\sqrt{\frac{\frac{1}{3}-\lambda}{1-\lambda}}$&$1$\\
 \hline
  $\mu=\frac{1}{2},\lambda\in\left[0,\frac{1}{3}\right),x=\displaystyle\frac{a+b}{2}-\frac{b-a}{2}\sqrt{\frac{\frac{1}{3}-\lambda}{1-\lambda}}$&$3$\\
 \hline
  $\displaystyle x=\frac{a+b}{2},\lambda=\frac{1}{3}$&$3$\\
  \hline
  \end{tabular}}$
  \end{tb}

In the following the estimations of $E(f;\lambda,\mu;x)$, involving derivatives of order greater than 1, are obtained for all special cases of parameters $\lambda,\mu\in[0,1]$ and $x\in\left[a,\frac{a+b}{2}\right]$ given in Table \ref{table1}:

\noindent{\bf Case 1.} If $\lambda=1$, then
\begin{align}& E(f;\lambda,\mu;x):=\displaystyle\frac{1}{b-a}\int_a^bf(t)dt-\frac{f(a)+f(b)}{2}=\int_a^bK(t)f^{\prime\prime}(t)dt,\label{z}\\
&K(t)=\displaystyle\frac{1}{2(b-a)}(b-t)(a-t).\nonumber\end{align}
It follows $|E(f;\lambda,\mu;x)|\leq C\cdot\|f^{\prime\prime}\|_{\infty}$, where $C=\displaystyle\frac{(b-a)^2}{12}$.

The quadrature formula (\ref{z}) is known as {\it trapezoidal rule} and the functional $E$ has also the following representation using the first derivative of $f$:
\begin{align*}& E(f;\lambda,\mu;x)=-\int_a^bK_1(t)f^{\prime}(t)dt, \textrm{ where }\\
&K_1(t)=\displaystyle\frac{1}{(b-a)}\left(t-\frac{a+b}{2}\right) \textrm{ and }\\
&|E(f;\lambda,\mu;x)|\leq \displaystyle\frac{b-a}{4}\|f^{\prime}\|_{\infty}.\end{align*}

\noindent{\bf Case 2.} If $x=\displaystyle\frac{a+b}{2}$ and $\lambda\ne \displaystyle\frac{1}{3}$, then
\begin{align}& E(f;\lambda,\mu;x):=\displaystyle\frac{1}{b-a}\int_a^bf(t)dt-\left\{\lambda\frac{f(a)+f(b)}{2}+(1-\lambda)f\left(\frac{a+b}{2}\right)\right\}
=\int_a^bK(t)f^{\prime\prime}(t)dt,\label{u}\\
&K(t)=\left\{\begin{array}{l}\displaystyle\frac{1}{2(b-a)}(t-a)\left[(t-a)-\lambda(b-a)\right],\, t\in\left[a,\frac{a+b}{2}\right]\nonumber,\\
\vspace{-0.3cm}\\
\displaystyle\frac{1}{2(b-a)}(b-t)\left[(b-t)-\lambda(b-a)\right],\, t\in\left(\frac{a+b}{2},b\right].\end{array}\right.\nonumber
\end{align}
We have $|E(f;\lambda,\mu;x)|\leq C\cdot\|f^{\prime\prime}\|_{\infty}$, where
 $C=\left\{\begin{array}{l}\displaystyle\frac{1}{24}(b-a)^2(8\lambda^2-3\lambda+1),\,\lambda\in\left[\displaystyle 0,\frac{1}{2}\right),\\
 \vspace{-0.3cm}\\
 \displaystyle\frac{1}{24}(b-a)^2(3\lambda-1),\,\lambda\in\left[\displaystyle \frac{1}{2},1\right].\end{array}\right.$

The functional (\ref{u}) has also the following representation using the  first derivative
\begin{align*}&E(f;\lambda,\mu;x):=\displaystyle\int_a^bK_1(t)f^{\prime}(t)dt,\textrm{ where} \\
&K_1(t)=\left\{\begin{array}{l} \displaystyle\frac{1}{2(b-a)}\left[2a-2t+\lambda(b-a)\right], t\in\left[a,\frac{a+b}{2}\right]\\
\vspace{-0.3cm}\\
 \displaystyle\frac{1}{2(b-a)}\left[2b-2t-\lambda(b-a)\right], t\in\left(\frac{a+b}{2},b\right]\\
\end{array}\right.\\
&|E(f;\lambda,\mu;x)|\leq\displaystyle\frac{b-a}{4}(2\lambda^2-2\lambda+1)\|f^{\prime}\|_{\infty}.
\end{align*}
We remark that for $\lambda=0$ we  obtain the well-known {\it mid-point rule}:
$$E(f;\lambda,\mu,x):=\displaystyle\frac{1}{b-a}\int_{a}^{b}f(x)dx-f\left(\frac{a+b}{2}\right)  $$
with the following representation of the remainder term:
\begin{itemize}
\item[a)] $E(f;\lambda,\mu;x):=\displaystyle\int_a^b K(t)f^{\prime\prime}(t)dt,$\\
$K(t):=\left\{\begin{array}{l}
\displaystyle\frac{(t-a)^2}{2(b-a)}, t\in\left[ a,\frac{a+b}{2}\right],\\
\vspace{-0.3cm}\\
\displaystyle\frac{(b-t)^2}{2(b-a)}, t\in\left(\frac{a+b}{2},b\right],\\
\end{array}\right.$\\
$|E(f;\lambda,\mu;x)|\leq\displaystyle\frac{(b-a)^2}{24}\| f^{\prime\prime}\|_{\infty}.$
\item[b)] $E(f;\lambda,\mu;x):=\displaystyle\int_a^b K_1(t)f^{\prime}(t)dt,$\\
$K_1(t):=\left\{\begin{array}{l}
\displaystyle\frac{a-t}{b-a}, t\in\left[ a,\frac{a+b}{2}\right],\\
\vspace{-0.3cm}\\
\displaystyle\frac{b-t}{b-a}, t\in\left(\frac{a+b}{2},b\right],\\
\end{array}\right.$\\
$|E(f;\lambda,\mu;x)|\leq\displaystyle\frac{b-a}{4}\| f^{\prime}\|_{\infty}.$
\end{itemize}

For $\lambda=\displaystyle\frac{1}{2}$ we obtain remainder term of the composite trapezoidal rule at $a,\displaystyle\frac{a+b}{2}$,
 and $b$, namely
$$ E(f;\lambda,\mu;x)=\displaystyle\frac{1}{b-a}\int_a^b f(x)dx-\left\{\displaystyle\frac{f(a)+f(b)}{4}+\frac{1}{2}f\left(\frac{a+b}{2}\right)\right\} $$
with
$$ \left| E(f;\lambda,\mu;x) \right|\leq \displaystyle\frac{1}{48}(b-a)^2\|f^{\prime\prime}\|_{\infty},f\in C^2[a,b] $$
and
$$|E(f;\lambda,\mu;x)|\leq \frac{1}{8}(b-a)\| f^{\prime}\|_{\infty},f\in C^1[a,b].$$

\noindent{\bf Case 3.} If $\mu=\displaystyle\frac{1}{2}$, $\lambda\in\left[0,\displaystyle\frac{1}{3}\right)$ and
$x\ne \displaystyle\frac{1}{2}\left[a+b-(b-a)\sqrt{\frac{\frac{1}{3}-\lambda}{1-\lambda}}\right]$, then
\begin{align}\label{E}& E(f;\lambda,\mu;x):=\displaystyle\frac{1}{b-a}\int_a^bf(t)dt-\left\{\lambda\frac{f(a)+f(b)}{2}+(1-\lambda)\cdot\frac{f(x)+f(a+b-x)}{2}\right\}
\nonumber\\ &=\int_a^bK(t)f^{\prime\prime}(t)dt,\\
&K(t)=\left\{\begin{array}{l}
\displaystyle\frac{(t-a)}{2(b-a)}\left[(t-a)-\lambda(b-a)\right],\,t\in\left[a,x\right],\\
\vspace{-0.3cm}\\
\displaystyle\frac{(b-t)(a-t)}{2(b-a)}-(1-\lambda)\frac{(a-x)}{2},\,t\in\left(x,a+b-x\right),\\
\vspace{-0.3cm}\\
\displaystyle\frac{(b-t)}{2(b-a)}\left[(b-t)-\lambda(b-a)\right],\,t\in\left[a+b-x,b\right].\\
\end{array}\right.\nonumber\end{align}
It follows $|E(f;\lambda,\mu;x)|\leq C\cdot\|f^{\prime\prime}\|_{\infty}$, where
\begin{equation}\label{ec}
C=\left\{\begin{array}{l}
\displaystyle\frac{(b-a)^2}{12}+\frac{1}{2}(1-\lambda)(a-x)(b-x),\,x\in[a,a+\lambda(b-a)],\\
\vspace{-0.3cm}\\
\displaystyle \frac{1}{6}(b-a)^2\left[1+\frac{4(1-\lambda)(a-x)}{b-a}\right]^{3/2}-\frac{1}{12}(b-a)^2(1-4\lambda^3)-\frac{1}{2}(1-\lambda)(a-x)(b-x),\\
\vspace{-0.3cm}\\
\displaystyle x\in\left(a+\lambda(b-a),a+\frac{b-a}{4(1-\lambda)}\right),\\
\vspace{-0.3cm}\\
\displaystyle\frac{(4\lambda^3-1)(b-a)^2}{12}+\frac{1}{2}(1-\lambda)(b-x)(x-a),\,x\in\left[a+\frac{b-a}{4(1-\lambda)},\frac{a+b}{2}\right].
\end{array}\right.
\end{equation}

\noindent{\bf Case 4.} If $\mu=\displaystyle\frac{1}{2}$ and $\lambda\in\left(\displaystyle\frac{1}{3},1\right)$, the functional $E(f;\lambda,\mu;x)$ is defined in
 (\ref{E}) and
 $|E(f;\lambda,\mu;x)|\leq C\cdot\|f^{\prime\prime}\|_{\infty}$, where for $\lambda\in\left(\displaystyle\frac{1}{3},\frac{1}{2}\right)$ the constant $C$ is defined in
 (\ref{ec}) and for $\lambda\in\left(\displaystyle\frac{1}{2},1\right)$ we have
$C=\displaystyle\displaystyle\frac{(b-a)^2}{12}+\frac{1}{2}(1-\lambda)(b-x)(a-x).$

\noindent{\bf Case 5.}  If $\mu=\displaystyle\frac{1}{2}$, $\lambda\in\left[0,\displaystyle\frac{1}{3}\right)$ and
$x= \displaystyle\frac{1}{2}\left[a+b-(b-a)\sqrt{\frac{\frac{1}{3}-\lambda}{1-\lambda}}\right]$, then
\begin{align*}& E(f;\lambda,\mu;x):=\displaystyle\frac{1}{b-a}\int_a^bf(t)dt-\left\{\lambda\frac{f(a)+f(b)}{2}
+(1-\lambda)\cdot\frac{f(x)+f(a+b-x)}{2}\right\}\\
&=\frac{1}{3!}\int_a^bK(t)f^{(4)}(t)dt,\\
&K(t)=\left\{\begin{array}{l}
\displaystyle\frac{(t-a)^3}{4(b-a)}\left[t-2\lambda(b-a)-a\right],\,t\in\left[a,x\right],\\
\vspace{-0.3cm}\\
\displaystyle\frac{1}{b-a}\cdot\frac{(b-t)^4}{4}-\left\{\lambda\frac{(b-t)^3}{2}+(1-\lambda)\cdot\frac{(a+b-x-t)^3}{2}\right\},\,
t\in\left(x,a+b-x\right),\\
\vspace{-0.3cm}\\
\displaystyle\frac{(b-t)^3}{4(b-a)}\left[b-2\lambda(b-a)-t\right],\,t\in\left[a+b-x,b\right].\\
\end{array}\right.\end{align*}
We have $|E(f;\lambda,\mu;x)|\leq C\cdot\|f^{(4)}\|_{\infty}$, where
$C=\displaystyle\frac{1}{720(1-\lambda)}\left(\lambda-\frac{1}{6}\right)\left(b-a\right)^4$, for $\displaystyle\lambda\in\left[\frac{1}{4},\frac{1}{3}\right)$.

For $\lambda=\displaystyle\frac{1}{4}$ we have the {\it Newton-Simpson's quadrature rule} ({\it Simpson's 3/8 rule}):
\begin{align}\label{r}
E(f;\lambda,\mu;x)&:=\displaystyle\frac{1}{b-a}\int_a^bf(t)dt-\frac{1}{8}\left\{ f(a)+f(b)+3f\left(\frac{2a+b}{3}\right)+3f\left(\frac{a+2b}{3}\right)\right\}\\
&=\frac{1}{3!}\int_a^bK(t)f^{(4)}(t)dt,\nonumber
\end{align}
where
$$ K(t)=\displaystyle\left\{\begin{array}{l}
\displaystyle\frac{(t-a)^3}{4(b-a)}\left(t-\frac{a+b}{2}\right),\, t\in\left[a, \frac{2a+b}{3}\right],\\
\vspace{-0.3cm}\\
\displaystyle \frac{1}{72(b-a)}\left(6t^2-6t(a+b)+a^2+b^2+4ab\right)\cdot \left(3t^2-3t(a+b)+a^2+ab+b^2 \right),\\
\vspace{-0.3cm}\\
\displaystyle t\in\left(\frac{2a+b}{3},\frac{a+2b}{3}\right),\\
\vspace{-0.3cm}\\
\displaystyle\frac{(b-t)^3}{4(b-a)}\left(\frac{a+b}{2}-t\right),\,t\in\left[\frac{a+2b}{3},b\right],
\end{array}
\right.$$
and $\left|E(f;\lambda,\mu;x)\right|\leq\displaystyle\frac{1}{6480}\|f^{(4)}\|_{\infty}$.

The functional (\ref{r}) has also the following representations:
\begin{itemize}
\item [a)] $E(f;\lambda,\mu;x):=-\displaystyle\frac{1}{3!}\int_{a}^bK_1(t)f^{(3)}(t)dt,\,\textrm{where}$\\
$K_1(t)=\left\{\begin{array}{l}
\displaystyle\frac{(t-a)^2}{b-a}\left(t-\frac{5a+3b}{8}\right),\, t\in\left[a,\frac{2a+b}{3}\right]\\
\vspace{-0.3cm}\\
\displaystyle\frac{1}{b-a}\left(t-\frac{a+b}{2}\right)^3, t\in\left(\frac{2a+b}{3},\frac{a+2b}{3}\right)\\
\vspace{-0.3cm}\\
\displaystyle\frac{(b-t)^2}{b-a}\left(t-\frac{3a+5b}{8}\right), t\in\left[\frac{a+2b}{3},b\right]\\
\end{array}\right.$\\
$\left|E(f;\lambda,\mu;x)\right|\leq \displaystyle\frac{1}{1728}(b-a)^3\|f^{(3)}\|_{\infty}$
\item[b)]$\left|E(f;\lambda,\mu;x)\right|= \displaystyle\frac{1}{3!}\int_a^bK_2(t)f^{(2)}(t)dt,\textrm{ where }$\\
$K_2(t)=\left\{\begin{array}{l}
\displaystyle\frac{3(t-a)}{b-a}\left(t-\frac{3a+b}{4}\right),\, t\in\left[a,\frac{2a+b}{3}\right],\\
\vspace{-0.3cm}\\
\displaystyle\frac{3}{b-a}\left(t-\frac{a+b}{2}\right)^2, t\in\left(\frac{2a+b}{3},\frac{a+2b}{3}\right),\\
\vspace{-0.3cm}\\
\displaystyle\frac{3(b-t)}{b-a}\left(\frac{a+3b}{4}-t\right), t\in\left[\frac{a+2b}{3},b\right],\\
\end{array}\right.$\\
$E(f;\lambda,\mu;x)|\leq\displaystyle\frac{1}{192}(b-a)^2\|f^{(2)}\|_{\infty}.$\\
\item[c)]$\left|E(f;\lambda,\mu;x)\right|= -\displaystyle\frac{1}{3!}\int_a^bK_3(t)f^{\prime}(t)dt,\textrm{ where }$\\
$K_3(t)=\left\{\begin{array}{l}
\displaystyle\frac{6}{b-a}\left(t-\frac{7a+b}{8}\right),\, t\in\left[a,\frac{2a+b}{3}\right],\\
\vspace{-0.3cm}\\
\displaystyle\frac{6}{b-a}\left(t-\frac{a+b}{2}\right)^2, t\in\left(\frac{2a+b}{3},\frac{a+2b}{3}\right),\\
\vspace{-0.3cm}\\
\displaystyle\frac{6}{b-a}\left(t-\frac{a+7b}{8}\right), t\in\left[\frac{a+2b}{3},b\right],\\
\end{array}\right.$\\
$E(f;\lambda,\mu;x)|\leq\displaystyle\frac{25}{288}(b-a)\|f^{\prime}\|_{\infty}.$
\end{itemize}

\noindent{\bf Case 6.} If $x=\displaystyle\frac{a+b}{2}$ and $\lambda= \displaystyle\frac{1}{3}$, then
\begin{align}& E(f;\lambda,\mu;x):=\displaystyle\frac{1}{b-a}\int_a^bf(t)dt-\left\{\frac{f(a)+f(b)}{6}+\frac{2}{3}f\left(\frac{a+b}{2}\right)\right\}
=\frac{1}{3!}\int_a^bK(t)f^{(4)}(t)dt,\label{vv}\\
&K(t)=\left\{\begin{array}{l}
\displaystyle\frac{1}{12(b-a)}(a-t)^3(a-3t+2b),\,t\in\left[a,\frac{a+b}{2}\right],\\
\vspace{-0.3cm}\\
\displaystyle\frac{1}{12(b-a)}(b-t)^3(b-3t+2a),\,t\in\left(\frac{a+b}{2},b\right].\\
\end{array}\right.\nonumber\end{align}
Then, we have $|E(f;\lambda,\mu;x)|\leq C\cdot\|f^{(4)}\|_{\infty}$, where $C=\displaystyle\frac{1}{2880}(b-a)^4$.

The quadrature formula in  (\ref{vv}) is known as {\it Simpson's rule}  and the functional $E$ has also the following representations:
\begin{itemize}
\item[a)] $E(f;\lambda,\mu;x)=\displaystyle\frac{1}{3!}\int_a^b K_1(t)f^{(3)}(t)dt,\textrm{ where}$\\
$ K_1(t)=\left\{\begin{array}{l}\displaystyle\frac{1}{2(b-a)}(a-t)^2(a+b-2t), t\in\left[ a,\frac{a+b}{2}\right],\\
\vspace{-0.3cm}\\
\displaystyle\frac{1}{2(b-a)}(b-t)^2(a+b-2t), t\in\left(\frac{a+b}{2},b\right],\end{array}\right. $\\
$\left|E(f;\lambda,\mu;x)\right|\leq \displaystyle\frac{1}{576}(b-a)^3\| f^{(3)}\|_{\infty}.$
\item[b)] $E(f;\lambda,\mu;x)=\displaystyle\frac{1}{3!}\int_a^b K_2(t)f^{(2)}(t)dt,\textrm{ where}$\\
$ K_2(t)=\left\{\begin{array}{l}\displaystyle\frac{a-t}{b-a}(2a+b-3t), t\in\left[ a,\frac{a+b}{2}\right],\\
\vspace{-0.3cm}\\
\displaystyle\frac{b-t}{b-a}(a+2b-3t), t\in\left(\frac{a+b}{2},b\right],\end{array}\right. $\\
$\left|E(f;\lambda,\mu;x)\right|\leq \displaystyle\frac{1}{81}(b-a)^2\| f^{\prime\prime}\|_{\infty}.$
\item[c)]  $E(f;\lambda,\mu;x)=\displaystyle\frac{1}{3!}\int_a^b K_3(t)f^{\prime}(t)dt,\textrm{ where}$\\
$ K_3(t)=\left\{\begin{array}{l}\displaystyle\frac{1}{b-a}(b-6t+5a), t\in\left[ a,\frac{a+b}{2}\right],\\
\vspace{-0.3cm}\\
\displaystyle\frac{1}{b-a}(a+5b-6t), t\in\left(\frac{a+b}{2},b\right],\end{array}\right. $\\
$\left|E(f;\lambda,\mu;x)\right|\leq \displaystyle\frac{5}{36}(b-a)\| f^{\prime}\|_{\infty}.$
\end{itemize}

\section{Estimates for all continuous functions}
The error functional $E(f;\lambda,\mu;x)$ is well-defined for all functions continuous on $[a,b]$. However, all previous considerations only dealt with differentiable ones.
We will now turn the previous inequalities into such for $f\in C[a,b]$ in general.

We first describe the $K$-functional technique in a way suited to our situation and later discuss several special cases which
will depend on the maximum degrees of exactness.

For a compact interval $[a,b]\subset \mathbb{R}$, $a<b$, we consider a bounded linear functional $L:C[a,b]\to \mathbb{R}$, $L\ne 0$.
We assume that for all functions $g$ in $C^s[a,b]$, $s\geq 1$, we have inequalities of the form
$$ |L(g)|\leq C_{L}\cdot \|g^{(s)}\|_{\infty}, $$
where $C_{L}$ is independent of $g$. If $f\in C[a,b]$ is fixed, then for all $g\in C^{s}[a,b]$ we write
\begin{align*}
|L(f)|&=|L(f-g+g|\leq |L(f-g)|+|L(g)|\\
&\leq \| L\|\cdot \| f-g\|_{\infty}+C_{L}\cdot \|g^{(s)}\|_{\infty}\\
&=\| L\|\cdot \left(\| f-g\|_{\infty}+\displaystyle\frac{C_{L}}{\| L\|}\cdot \|g^{(s)}\|_{\infty}\right).
\end{align*}
This implies
\begin{align}
|L(f)|&\leq \|L\|\cdot\inf\left\{\|f-g\|_{\infty}+\displaystyle\frac{C_{L}}{\| L\|}\cdot\| g^{(s)}\|_{\infty}: g\in C^{s}[a,b]\right\}\label{F1}\\
&=:\|L\|\cdot K\left(f,\displaystyle\frac{C_{L}}{\| L\|}; C[a,b], C^s[a,b]\right).\nonumber
\end{align}
The latter is the $K$-functional with respect to the spaces $(C[a,b],\|\cdot\|_{\infty})$
and $(C^s[a,b],\|\cdot^{(s)}\|_{\infty})$.

For the error functional $E(f;\lambda;\mu;x)$ of the generalized Alomari functional it is immediately seen that
$$ \| E(\cdot;\lambda;\mu;x)\|\leq 2. $$
Combining this with (\ref{e6}) in Theorem \ref{t5} we arrive at
$$| E(f;\lambda,\mu;x)\|\leq 2\cdot K\left(f,\displaystyle\frac{C(x,\lambda,\mu)}{4(b-a)};C, C^1\right)  $$
for all $f\in C[a,b],\lambda,\mu\in [0,1]$ and $x\in \left[a,\displaystyle\frac{a+b}{2}\right]$.
 An analogous inequality holds for $x\in \left(\displaystyle\frac{a+b}{2},b\right]$.

 \subsection{Estimates in terms of the least concave majorant of the first order modulus of continuity}
 The following description is taken from \cite{Pal2010}.

 Recall that the least concave majorant of a function $f:[a,b]\to \mathbb{R}$ is the function $\tilde{f}:[a,b]\to \mathbb{R}$, defined by
 $$ \tilde{f}(x)=\inf\left\{g(x),g\geq f, g\textrm{ concave }\right\}. $$
 From this definition it follows that
 $$ \tilde{f}(x)=\inf\{l(x),l\geq f, l\textrm{ linear }\} $$
 and
 $$ \tilde{f}(x)=\sup\left\{\displaystyle\frac{(t-x)f(s)+(x-s)f(t)}{t-s},a\leq s\leq x\leq t\leq b, s<t\right\}. $$
 The first modulus of continuity of a bounded function $f:[a,b]\to \mathbb{R}$ is defined by
 $$ \omega(f,h):=\sup\left\{|f(u)-f(v)|,u,v\in [a,b],|u-v|\leq h\right\},\textrm{ for } h>0. $$
 If, for a given function $f$ we construct the least concave majorant of the map $t\to \omega(f,t)$, $t\in [0,b-a]$,
  we obtain the least concave majorant of the first modulus of $f$, which is denote by $\tilde{\omega}(f,\cdot)$. Then we have
  \begin{thm}For any $f\in C[a,b]$ and any $0<t\leq (b-a)/2$,
  $$K(f,t; C[a,b], C^1[a,b])=\displaystyle\frac{1}{2}\tilde{\omega}(f,2t).  $$
  \end{thm}
  Now recall (\ref{F1}). This leads immediately to
  \begin{align*}
  \left| E(f,\lambda,\mu;x)\right|\leq 2\cdot K\left(f,\displaystyle\frac{C(x,\lambda,\mu)}{2};C[a,b], C^1[a,b] \right) 
  =\tilde{\omega}\left(f,\displaystyle {C(x,\lambda,\mu)}\right),
  \end{align*}
  for all $\lambda,\mu\in[0,1], x\in\left[a,\displaystyle\frac{a+b}{2}\right]$, $f\in C[a,b]$.
  Note that $\displaystyle\frac{C(x,\lambda,\mu)}{2}\leq\frac{b-a}{2}$.

As an example we reconsider the midpoint rule in which we had
$$ \left|E\left(f;0,\mu,\displaystyle\frac{a+b}{2}\right)\right|\leq\displaystyle\frac{b-a}{4}\| f^{\prime}\|_{\infty},\, f\in C^1[a,b]. $$
We thus see that for all $f\in C[a,b]$ we obtain
$$ \left| E\left(f;0,\mu,\frac{a+b}{2}\right)\right|\leq 2K\left(f,\frac{b-a}{8};C[a,b],C^1[a,b]\right)={\tilde\omega}\left(f,\frac{b-a}{4}\right), $$
which exactely implies what we invested for $C^1$ functions.

   \subsection{Estimates in terms of the second modulus of smoothness}
   In the above Cases 1.-4. we gave estimates of the type
   $$ | E(g,\lambda,\mu,x)|\leq C\cdot\|g^{\prime\prime}\|_{\infty},\textrm{ for } g\in C^2[a,b], $$
   and for special values of $\lambda, \mu$ and $x$. We will now show how these imply estimates for arbitrary $f\in C[a,b]$.

   Recalling that for such $f$ the second modulus is given by
   $$ \omega_2(f;h)=\sup\left\{|f(x-t)-2f(x)+f(x+t)|:x\pm t\in[a,b],h\geq t\right\}, $$
   it was shown in \cite{GoKo1994} (see Lemmas 2.1 and 4.1 there) that for $0< t\leq \displaystyle\frac{1}{2}(b-a)$ one has
   $$ K\left(f,t;C[a,b], C^2[a,b]\right)\leq \displaystyle\frac{9}{4}\cdot\omega_2\left(f;\sqrt{t}\right). $$
   Combining this with (\ref{F1}) we obtain in cases 1.-4. the inequalities
   $$\left|E(f,\lambda,\mu,x)\right|\leq 4.5\cdot \omega_{2}\left(f;\sqrt{\displaystyle\frac{C}{2}}\right),\,f\in C[a,b], $$
   where the constants $C$ are those from the four cases (and subcases).
   In order to fully justify this statement, one has to show that always
   $$ \sqrt{\displaystyle\frac{C}{2}}\leq \displaystyle\frac{b-a}{2}, $$
   which can be done by inspection.

Again we consider the midpoint rule. There it was show that
$$   \left| E\left(f;0,\mu,\frac{a+b}{2}\right) \right|\leq \frac{(b-a)^2}{24}\| f^{\prime\prime}\|_{\infty}, f\in C^2[a,b].$$
So from the above we get
$$ \left| E\left(f;0,\mu,\frac{a+b}{2}\right)\right|\leq 4.5\,\omega_2\left(f;\displaystyle\frac{b-a}{4\sqrt{3}}\right). $$
The latter inequality implies such in terms of $\| f^{\prime}\|_{\infty}$ and $\| f^{\prime\prime}\|_{\infty}$, but-as far as the numerical constants are  concerned- these are worse than what we already knew for $C^1$ and $C^2$, respectively.

\subsection{Estimates in terms of the fourth order modulus of smoothness}

In Cases 5. and 6. we have a degree of exactness of three, that is, inequalities of the form
$$ |E(g;\lambda,\mu;x)|\leq C \cdot \| g^{(4)}\|_{\infty}, g\in C^{4}[a,b].  $$

Now (\ref{F1}) implies
$$ | E(f;\lambda,\mu;x)| \leq 2\cdot K\left(f,\displaystyle\frac{C}{2};C[a,b],C^4[a,b]\right).$$
This can be used to obtain inequalities in terms of the fourth order modulus, but-as will be seen below- the situation
 is becoming more and more inconvenient due to the growth of the constants (known to us). Also, we are only aware of a result for $[a,b]=[-1,1]$.
 The following theorem is taken from the diploma thesis of Angelika Sperling \cite{Spe1984}.
 \begin{thm} (see \cite{Spe1984}, Bemerkung 10.22 (ii)) For  $f\in C[-1,1]$, $s\in \mathbb{N}$ and $0<t\leq \displaystyle\frac{2}{s}$ one has
 $$ K(f,t^s;C[-1,1], C^s[-1,1])\leq (2s)^{2s}\cdot D(s)\cdot\omega_{s}(f;t), $$ where
 \begin{align*}
 &D(s)=2s^s+\tilde{D}_{1}(s)+\displaystyle\frac{1}{2^ss^{2s}}\tilde{D}_{2}(s)\\
 & \tilde{D}_{1}(s)=(2^s-1)(2\cdot D_1(s)+3)\\
 &\tilde{D}_{2}(s)=2s^s\left(D_2(s)+{s-1\choose \left[\frac{s-1}{2}\right]}\cdot\frac{(2s+1)!}{(s+1)!}\right)\\
 &D_1(s)=\displaystyle\sum_{k=1}^{s-1}{s\choose k}\cdot\displaystyle\frac{(2s+1)!}{(2s+1-k)!}\cdot{k-1\choose\left[\frac{k-1}{2}\right]}, s\geq 2\\
& D_2(s)=\displaystyle\sum_{k=1}^{s-1}{s\choose k}\cdot\frac{(2s+1)!}{(2s+1-k)!}{k-1\choose \left[\frac{k-1}{2}\right]}\cdot m_{s,s-k},s\geq 2,\\
&m_{s,i}=\left\{\begin{array}{l} (2^s-1), i=0,\\
(2^s-1)(s+1)^i\displaystyle\prod_{j=1}^i\left(1+\ln(s-j\right)),1\leq i\leq s-1, 0\leq i\leq s-1.\end{array}\right.
 \end{align*}
 \end{thm}
For the case $s=4$ it was show in Bemerkung 10.11 (iii) of \cite[p.~146]{Spe1984}  that
\begin{align*}
{\tilde D}_1(4)=135\, 045\textrm{ and }
{\tilde D}_{2}(4)=847\,078\,494.
\end{align*}
Hence
$$ D(4)=2\cdot 4^4+135\,045+\displaystyle\frac{1}{2^4\cdot 4^8}\cdot 847\, 078\, 494< 136\, 365. $$
The suitable example is now Case 6 from above (Simpson's rule). Combining the previous information with Speling's inequality yields
\begin{align*}E\left(f;\displaystyle\frac{1}{3},\mu,\frac{a+b}{2}\right)&\leq 2K\left(f,\frac{1}{2}\cdot\frac{1}{2880}\cdot 2^4; C^0, C^4\right)\\
&\leq 2\cdot(2\cdot 4)^8\cdot D(4)\cdot\omega_4\left(f;\sqrt[4]{\frac{1}{2}\cdot\frac{1}{2880}\cdot 2^4}\right)\\
&=4\,575\,650\,119\,680\cdot\omega_4(f;0.2295\dots).
\end{align*}

\noindent{\bf Acknowledgment.} Project financed from Lucian Blaga
University of Sibiu research grants LBUS-IRG-2015-01, No.2032/7. It
was carried out while the second author was visiting LBUS under an
Erasmus + grant.

$  $

\end{document}